\documentclass[10pt,a4paper]{article}
\usepackage[english]{babel}
\usepackage{polski}
\usepackage{amsmath, amsthm, amsfonts}
\usepackage{tikz-cd}
\usepackage{amssymb}

\usepackage{stackengine}
\usepackage{scalerel}

\theoremstyle{plain}
\newtheorem{theorem}{Theorem}[section]
\newtheorem{lemma}[theorem]{Lemma}
\newtheorem{proposition}[theorem]{Proposition}
\theoremstyle{definition}
\theoremstyle{definition}
\newtheorem{remark}[theorem]{Remark}
\theoremstyle{definition}
\theoremstyle{definition}
\newtheorem{example}[theorem]{Example}

\theoremstyle{definition}

\newtheorem{construction}[theorem]{Construction}
\usepackage[]{graphicx}

\title{Cox Rings and algebraic maps}

\begin{document}

\author{Tomasz Ma\'{n}dziuk\thanks{E-mail:~\textsf{t.mandziuk@mimuw.edu.pl}, Faculty of Mathematics, Computer Science and Mechanics, University of Warsaw, ul. Banacha 2, 02-097 Warszawa, Poland}}

\maketitle{}

\begin{abstract}
\noindent Given a morphism $F: X\rightarrow Y$ from a Mori Dream Space $X$ to a smooth Mori Dream Space $Y$ and quasicoherent sheaves $\mathcal{F}$ on $X$ and $\mathcal{G}$ on $Y$, we describe the inverse image of $\mathcal{G}$ by $F$ and the direct image of $\mathcal{F}$ by $F$ in terms of the corresponding modules over the Cox rings graded in the class groups.
\end{abstract}

\tableofcontents

\section*{Introduction}

In this paper, we will be interested in quasicoherent sheaves on Mori Dream Spaces. In \cite{Cox95b}, Cox described the so called \textbf{homogeneous coordinate ring} of a toric variety $X$. It is a $\mathrm{Cl}(X)$-graded ring, where $\mathrm{Cl}(X)$ is the divisor class group of $X$. The construction was later generalized to more general varieties and is now known as the \textbf{Cox ring}. It generalizes the homogeneous coordinate ring of a projective variety but it does not depend, up to a non-canonical isomorphism, on the choice of embedding in any affine or projective variety. For normal toric varieties, the Cox ring is always a $\mathrm{Cl}(X)$-graded polynomial ring. In general, a variety admitting a finitely generated Cox ring is called a \textbf{Mori Dream Space} (MDS). 

There is a correspondence between quasicoherent sheaves on a MDS $X$ and $\mathrm{Cl}(X)$-graded modules over the Cox ring of $X$. Suppose $X$ and $Y$ are MDSes with Cox rings $R$ and $S$, respectively, and $\mathcal{F}$ is a quasicoherent sheaf on $X$, $\mathcal{G}$ is a quasicoherent sheaf on $Y$. Assume $F:X\rightarrow Y$ is a morphism. Let $\mathcal{F}$ correspond to a $\mathrm{Cl}(X)$-graded $R$-module $M$ and let $\mathcal{G}$ correspond to a $\mathrm{Cl}(Y)$-graded $S$-module $N$. One may ask, to which $\mathrm{Cl}(X)$-graded $R$-module corresponds the inverse image sheaf $F^*\mathcal{G}$ and to which $\mathrm{Cl}(Y)$-graded $S$-module corresponds the direct image sheaf $F_*\mathcal{F}$. We answer these two questions in Theorems \ref{t:pullback} and \ref{t:pushforward}, respectively.

Since the Cox ring $R$ is $\mathrm{Cl}(X)$-graded, $\operatorname{Spec}R$ comes with an action of a quasitorus $H_X = \operatorname{Spec}(k[\mathrm{Cl}(X)])$, where $k$ is a fixed algebraically closed base field of characteristic zero. In \cite{Cox95b}, Cox proved also, that every normal toric variety $X$ can be obtained as the good quotient for the action of $H_X$ of an invariant open subset of $\operatorname{Spec}R$. Analogous result holds true for any MDS. Under some additional assumptions, every morphism of toric varieties can be lifted to affine varieties associated with their coordinate rings \cite{Cox95a}. We will use a similar result for MDSes.

In the first section we recall the relevant definitions and theorems. We start with the definition of Cox rings. Then we describe the quotient construction of MDSes. Subsequently, we introduce an affine open cover of MDSes that we will use for local study of quasicoherent sheaves. In the next two subsections we recall the results that are a basis of this paper. Namely, the correspondence between $\mathrm{Cl}(X)$-graded $R$-modules and quasicoherent sheaves on $X$ and the existence of a lift of a morphism of MDSes to a morphism of their Cox rings. In the second section we present the proofs of the main results: Theorem \ref{t:pullback} and Theorem \ref{t:pushforward} and give a few lemmas that will be used in the last section. The paper ends with two examples.

\subsection*{Acknowledgments}
This paper is based on my master thesis. I would like to thank my advisor, Jaros\l{}aw Buczy\'{n}ski, for introducing me to this topic and many helpful discussions that broaden my understanding of the related material. While writing the thesis I was supported by the Polish National Science Center (NCN), projects 2013/11/D/ST1/02580 and 2013/08/A/ST1/00804.

\section{Cox rings and Mori Dream Spaces}\label{r:pojecia}

All varieties that we will consider will be over a fixed algebraically closed field $k$ of characteristic zero. All constructions and definitions in this section are from \cite{ADHL15}. We will skip most of the proofs.

\subsection{Cox rings}
Let $X$ be a normal variety over $k$ with finitely generated class group $\mathrm{Cl}(X)$. We will recall the definition of the Cox sheaf $\mathcal{R}$ on $X$ and the Cox ring $R$ of X in two different settings using in each of them a 
different additional assumption. The idea is to construct a $\mathrm{Cl}(X)$-graded $\mathcal{O}_X$-algebra $\mathcal{R}$ such that $\mathcal{R}_{[D]}\cong \mathcal{O}_X(D)$. Then the Cox ring will be 
defined as $R=\bigoplus_{[D]\in \mathrm{Cl}(X)} \Gamma(X,\mathcal{R}_{[D]})$. The technical problem is the definition of the multiplication of sections of $\mathcal{R}$. 
We will first assume additionally that $\mathrm{Cl}(X)$ has no torsion.

\begin{construction}[Construction of the Cox ring. Version 1]
Let $X$ be a normal variety over $k$ with free finitely generated class group. Pick any subgroup $K\subset \mathrm{WDiv}(X)$ such that the quotient map $\mathrm{WDiv}(X)  \rightarrow \mathrm{Cl}(X)$ induces an isomorphism $c: K \rightarrow$ $\mathrm{Cl}(X)$. We define the \textbf{Cox sheaf} $\mathcal{R}$ as the quasicoherent sheaf of $\mathcal{O}_X$-algebras: $\mathcal{R} = \bigoplus_{E \in K} \mathcal{O}_X(E)$.
The structure of an $\mathcal{O}_X$-algebra on $\mathcal{R}$ comes from  multiplying the homogeneous sections in the function field $k(X)$. 
The  \textbf{Cox ring} is defined as: $R= \Gamma(X,\mathcal{R}) = \bigoplus_{E \in K} \Gamma(X, \mathcal{O}_{X}(E))$.
Up to isomorphism $\mathcal{R}$ does not depend on the choice of the subgroup $K$, see Construction 1.4.1.1 in \cite{ADHL15}.
\end{construction}
We will now present two easy examples of this construction.

\begin{example} Let $X$ be a normal affine variety with trivial class group. Then $\mathcal{R} = \mathcal{O}_X$ and the Cox ring of $X$ is the same as the affine coordinate ring of $X$.
\end{example}

\begin{example} Let $X = \mathbb{P}_{\mathbb{C}}^n$. 
Then $\mathrm{Cl}(X) \cong \mathbb{Z}H$ where $H$ is any hyperplane in 
$\mathbb{P}_{\mathbb{C}}^n$. 
We have $\Gamma(\mathbb{P}_{\mathbb{C}}^n,\mathcal{O}_{\mathbb{P}_{\mathbb{C}}^n}(dH)) \cong \mathbb{C}[x_0, ..., x_n]_d$ - the group of homogeneous polynomials of degree $d$ and these isomorphisms combine to give an isomorphism of the Cox ring of $\mathbb{P}_{\mathbb{C}}^n$ with the homogeneous coordinate ring of $\mathbb{P}_{\mathbb{C}}^n$. 
\end{example}

Note that the isomorphism class of the Cox ring of a variety  $X$ does not depend on the embedding of $X$ in any ambient variety. This is not the case for the homogeneous coordinate ring of a projective variety.

Requiring no torsion in $\mathrm{Cl}(X)$ is too restrictive. We will now remove this assumption, requiring instead that there are no non-constant global invertible functions on $X$, i.e. $\Gamma(X, \mathcal{O}_X^*) = k^*.$ This assumption is needed in the proof that the Cox sheaf (up to isomorphism) does not depend on the choices made in the following construction. Moreover, this additional assumption is easily satisfied, for instance if $X$ is projective or complete.

\begin{construction}[Construction of the Cox ring. Version 2]\label{k:pierscien}
Let $X$ be a normal variety over $k$ with finitely generated class group and $\Gamma(X,\mathcal{O}_X^*)=k^*$. Take any subgroup of the Weil divisor group $K\subset \mathrm{WDiv}(X)$ projecting onto $\mathrm{Cl}(X)$ under the quotient map $\mathrm{WDiv}(X) \rightarrow \mathrm{Cl}(X)$. Let $K^0$ be the kernel of $c:K\rightarrow \mathrm{Cl}(X)$ and let $\chi: K^0 \rightarrow k(X)^*$ be a character of $K^0$ such that ${\rm{div}}(\chi(E))=E$, for every $E$ in $K^0.$ Let $\mathcal{S}$ be the sheaf of $\mathcal{O}_X$-algebras associated with $K$, i.e. $\mathcal{S}=\bigoplus_{E \in K} \mathcal{O}_X(E)$.
Let $\mathcal{I}$ be the sheaf of ideals of $\mathcal{S}$ locally generated by the sections $1 - \chi(E)$ where $E \in K^0$. Here $1$ is a homogeneous element of degree $0$ and $\chi(E)$ is a homogeneous element of degree $-E$. Let $\pi : \mathcal{S} \rightarrow \mathcal{S}/\mathcal{I}$ be the projection map. The \textbf{Cox sheaf} is the quotient sheaf $\mathcal{R} = \mathcal{S}/\mathcal{I}$ with the $\mathrm{Cl}(X)$-grading given by:
\[ 
\mathcal{R} = \bigoplus\limits_{[D] \in \mathrm{Cl}(X)} \mathcal{R}_{[D]},\ \textrm{where } \mathcal{R}_{[D]} = \pi(\bigoplus\limits_{E \in c^{-1}([D])}\mathcal{S}_{E}).
\]
The \textbf{Cox ring} of $X$ is then given by: 
$R=\Gamma(X,\mathcal{R}) = \bigoplus_{[D]\in \mathrm{Cl}(X)}\Gamma(X, \mathcal{R}_{[D]})$.
Again, the Cox sheaf does not depend, up to isomorphism, on the choices of $K$ and $\chi$, see Proposition 1.4.2.2 in \cite{ADHL15}.
\end{construction}

In Lemma 1.4.3.4 in \cite{ADHL15} it is proved that for every $D \in K$, $\pi|_{\mathcal{S_D}} : \mathcal{S_D} \rightarrow \mathcal{R}_{[D]}$ is an isomorphism. Hence the Cox sheaf in either of the two constructions can be informally thought of as the direct sum of sheaves of $\mathcal{O}_X$-modules associated with each divisor class in $\mathrm{Cl}(X)$ with an appropriate $\mathcal{O}_X$-algebra structure. 

From now on, we restrict ourselves to considering the Cox rings of varieties fitting into the setting of the second construction. A normal variety $X$ with $\Gamma(X, \mathcal{O}_X^*)=k^*$ and a finitely generated class group and Cox ring will be called a \textbf{Mori Dream Space} (MDS). Note that this definition is not standard. For instance in \cite{ADHL15} it is assumed also that $X$ is projective but we do not need this assumption here.

\subsection{The quotient construction}
Every MDS can be constructed as the good quotient of an open invariant subset of an affine variety by a quasitorus action. In this section we will recall this construction following Construction 1.6.3.1 in \cite{ADHL15}.

\begin{construction}\label{k:iloraz}
Let $X$ be a normal variety over $k$ with finitely generated class group $\mathrm{Cl}(X)$ and no non-constant global invertible functions. Assume that the Cox ring $R$ is finitely generated. From Theorem 1.5.1.1 in \cite{ADHL15} it follows that $\mathcal{R}$ is a sheaf of reduced $\mathcal{O}_X$-algebras and from Proposition 1.6.1.1 in \cite{ADHL15} it follows that it is locally of finite type. Hence the relative spectrum of the Cox sheaf $\textbf{Spec}(\mathcal{R})$ is a variety. We will denote it by $\widehat{X}$ and call it the \textbf{characteristic space} of $X$. It comes with an action of the quasitorus $H_X$ associated with $\mathrm{Cl}(X)$ (i.e. $H_X = \operatorname{Spec}(k[\mathrm{Cl}(X)])$) and a good quotient for this action $\pi_X : \widehat{X} \rightarrow X$. Let $\overline{X}$ be the spectrum of the Cox ring. Since $R$ is integral and finitely generated as a $k$-algebra, $\overline{X}$ is a variety. We will call it the \textbf{total coordinate space} of $X$. Since $R$ is $\mathrm{Cl}(X)$-graded, the total coordinate space comes with an action of $H_X$. There is an equivariant open embedding $i_X : \widehat{X} \rightarrow \overline{X}$ with the complement of the image of codimension at least two. The homogeneous ideal of $R$ defining the complement $\overline{X}\setminus \widehat{X}$ will be denoted by $\mathcal{J}_{irr}(X)$ and will be called the \textbf{irrelevant ideal} of $X$.
\end{construction}

\subsection{The [D]-divisor and the [D]-localization}
In the study of local behaviour of quasicoherent sheaves on MDSes we will use the notions of a $[D]$-divisor and a $[D]$-localization from Section 1.5.2 in \cite{ADHL15}. In the notation from Construction \ref{k:pierscien}, take any divisor $D\in K$ and a non-zero $f \in \Gamma(X, \mathcal{R}_{[D]})$. 
Then by Lemma 1.4.3.3 in \cite{ADHL15} there exists a unique $\tilde{f}\in \Gamma(X, \mathcal{S}_{D})$ such that $\pi(\tilde{f})=f$. We define the \textbf{$[$D$]$-divisor} of $f$ as ${\rm{div}}_{[D]}(f) = {\rm{div}}(\tilde{f}) + D$. 
Note that this divisor is always effective. The $[D]$-divisor does not depend on the choice of a representative $D \in K$ and the choices made in Construction \ref{k:pierscien}. It follows easily from the definition that for $0 \neq f \in \Gamma(X, \mathcal{R}_{[D_1]})$ and $0 \neq g \in \Gamma(X, \mathcal{R}_{[D_2]})$ we have: 
\begin{equation}\label{r:baza}
{\rm{div}}_{[D_1]+[D_2]}(fg) = {\rm{div}}_{[D_1]}(f) + {\rm{div}}_{[D_2]}(g).
\end{equation}
For $0\neq f \in \Gamma(X, \mathcal{R}_{[D]})$ we define the \textbf{$[$D$]$-localization} of $X$ by $f$ as the complement of the support of the $[$D$]$-divisor of $f$, that is: $X_{[D],f} = X \setminus {\rm{Supp}}({\rm{div}}_{[D]}(f))$.
We will later need the following lemma.

\begin{lemma}\label{l:baza_sklejenie_2} Suppose $X$ is a MDS with the Cox ring $R$. Then for all divisor classes $[D],[E] \in \mathrm{Cl}(X)$ and for all non-zero $f \in R_{[D]}$ and $g \in R_{[E]}$ we have $X_{[D], f} \cap X_{[E], g} = X_{[D]+[E], fg}$.
\end{lemma}
\begin{proof}
Since both ${\rm{div}}_{[D]}f$ and ${\rm{div}}_{[E]}g$ are effective, equation (\ref{r:baza}) implies that 
\[
{\rm{Supp}}({\rm{div}}_{[D]+[E]}(fg))={\rm{Supp}}({\rm{div}}_{[D]}(f)) \cup {\rm{Supp}}({\rm{div}}_{[E]}(g)).
\]
\end{proof}

For a non-zero homogeneous element $f$ in $R$ we will denote by $R_{(f)}$ the degree zero part of $R_f$. 
Observe that if $X_{[D],f}$ is affine then, by Proposition 1.6.3.3 in \cite{ADHL15}, it is a good quotient of $\widehat{X}_f = \overline{X}_f \cong \operatorname{Spec}R_f$ by the action of $H_X$. In particular it is isomorphic to $\operatorname{Spec}R_{(f)}$.

We will later use the following lemma.

\begin{lemma}\label{l:baza_sklejanie} Suppose $X$ is a MDS with Cox ring $R$. Then the affine sets of the form $X_{[D], f}$ with $[D] \in \mathrm{Cl}(X)$ and $0 \neq f \in R_{[D]}$ form a basis for the topology of $X$. 
\end{lemma}
\begin{proof} Suppose that $U\subsetneq X$ is a non-empty open affine subset of $X$. We claim that the complement $X\setminus U$ of $U$ in $X$ is of pure codimension one. Suppose that $X\setminus U$ has an irreducible component $Z$ of codimension at least two. Let $V$ be an affine open subset of $X$ that intersects $Z$ but has empty intersection with all other irreducible components of $X\setminus U$. Since $X$ is separated, $U\cap V$ is affine. Thus, it is enough to prove the claim for affine $X$ and non-empty affine $U$ such that $X\setminus U$ is irreducible. We will denote it by $Z$. The inclusion $U=X\setminus Z \rightarrow X$ induces restriction morphism $\alpha:\Gamma(X, \mathcal{O}_X) \rightarrow \Gamma(X\setminus Z, \mathcal{O}_X)$. Since $X$ is integral, $\alpha$ is injective. By assumptions $Z$ has no points of codimension one. Thus, by Theorem 4.0.14 in \cite{CLS11},  $\alpha$ is also surjective. This gives a contradiction since both $U=X\setminus Z$ and $X$ are affine and the inclusion $U \rightarrow X$ is not an isomorphism. 

Since for every non-empty open affine subset $U \subsetneq X$ the complement of $U$ is of pure codimension 1, the lemma follows from Proposition 1.5.2.2 in \cite{ADHL15}
\end{proof}

\subsection{Quasicoherent sheaves on Mori Dream Spaces}
As in the case of toric varieties, there is a correspondence between quasicoherent sheaves on MDSes and modules over their Cox rings graded in the class group. Moreover, coherent sheaves correspond to the finitely generated modules.

\begin{proposition}[\cite{ADHL15}, 4.2.1.11]
Let $X$ be a Mori Dream Space with the Cox ring $R$. There is a functor:
\begin{center}
 \{$\mathrm{Cl}(X)$-graded $R$-modules\} $\rightarrow \normalfont{QCoh_X}$ given by
 $M \mapsto ({\pi_X}_*{i_X}^* \overline{M})_0$,
\end{center}
where $\overline{M}$ is the quasicoherent $\mathcal{O}_{\overline{X}}$-module associated with the $R$-module $M$. This functor is exact and essentially surjective. Moreover, it induces an exact and essentially surjective functor: 
 \begin{center}
 \{finitely generated $\mathrm{Cl}(X)$-graded $R$-modules\} $\rightarrow \normalfont{Coh_X}$.
 \end{center}
 \label{s:funktor}
\end{proposition}

By $\widetilde{M}$ we will denote the quasicoherent sheaf on $X$ corresponding to the $\mathrm{Cl}(X)$-graded $R$-module $M$ via the functor from the above proposition. A $\mathrm{Cl}(X)$-graded $R$-module $M$ also defines a quasicoherent sheaf on the total coordinate space $\overline{X}$. To make it clear which sheaf we are considering, we will adopt a non-standard convention of calling the latter $\overline{M}$.

We collect for further reference two facts that follow immediately from the proof of the above proposition.
 
\begin{proposition}\label{s:lematy_do_cofania}
Let $X$ be a Mori Dream Space with the Cox sheaf $\mathcal{R}$ and the Cox ring $R$. Let $\mathcal{F}$ be a quasicoherent sheaf on $X$. Denote by $M$, the $\mathrm{Cl}(X)$-graded $R$-module $\Gamma(\widehat{X},{i_X}^*\mathcal{F})$. Then the following statements hold true:
\begin{itemize}
\item[(i)] ${\pi_X}^*\widetilde{M} \cong {i_X}^* \overline{M}$,
\item[(ii)] ${\pi_X}_*{\pi_X}^* \mathcal{F} \cong \mathcal{F}\otimes_{\mathcal{O}_X}\mathcal{R}.$
\end{itemize}
\end{proposition}

\subsection{Lifting morphisms of Mori Dream Spaces to their Cox rings}
The main tool that will be used in the proofs of Theorems \ref{t:pullback} and \ref{t:pushforward} is the following result from \cite{HM16}.

\begin{theorem}\label{t:podnoszenie}
Let $X$ and $Y$ be Mori Dream Spaces with the Cox rings $R$ and $S$, respectively. Assume that $Y$ is smooth. Let $F:X\rightarrow Y$ be a morphism. Then there exists a morphism $\overline{F} : \overline{X} \rightarrow \overline{Y}$ such that:
\begin{itemize}
\item[(1)] the induced map on coordinate rings $\overline{F}^{*} : S \rightarrow R$ is a graded homomorphism of graded rings with respect to the pullback map $\mathrm{Cl}(Y) = \mathrm{Pic}(Y)  \rightarrow  \mathrm{Pic}(X) \rightarrow \mathrm{Cl}(X)$, and
\item[(2)] the following diagram is commutative:
\end{itemize}
\begin{center}
\begin{tikzcd}[row sep=normal, column sep=normal]
\overline{X} \arrow[dotted]{r}{\overline{F}} & \overline{Y} \\
\widehat{X} \arrow[hook]{u}{i_X} \arrow[dotted]{r}{\widehat{F}} \arrow[two heads]{d}{\pi_X} &
\widehat{Y} \arrow[two heads]{d}{\pi_Y} \arrow[hook]{u}{i_Y} \\
X \arrow{r}{F} 	& Y
\end{tikzcd}
\end{center}
\noindent where $\widehat{F}$ is the restriction of $\overline{F}$ to $\widehat{X}$.
\end{theorem}

\begin{remark} In the proofs of Theorems \ref{t:pullback} and \ref{t:pushforward} we will not explicitly use the assumption that $Y$ is smooth. We require only the existence of a lifting giving a commutative diagram as in Theorem \ref{t:podnoszenie} (2).
\end{remark}

\begin{remark} There are similar results on existence of a lift of a map of MDSes to Cox rings. In \cite{Cox95a} for a morphism from a complete toric variety into a smooth toric variety without torus factors. In \cite{BB13} there are considered rational maps of toric varieties. In this case the lift to the total coordinate spaces is a multi-valued function. It is generalized further in \cite{BK16} by considering rational maps of MDSes. See also \cite{HM16}.
\end{remark}

\begin{remark}
In general, for a morphism of MDSes $X \rightarrow Y$, there does not exist a lift to the total coordinate spaces as in Theorem \ref{t:podnoszenie} (2). A simple example is given in section 1.1.2 in \cite{BB13}.
\end{remark}

\section{Main results}\label{s_3}
In this section, we are in the following setup. Let $X$, $Y$ be Mori Dream Spaces. The Cox sheaves of $X$ and $Y$ are $\mathcal{R}$ and $\mathcal{S}$, respectively. The Cox rings of $X$ and $Y$ are $R$ and $S$, respectively. We have a morphism $F:X\rightarrow Y$. We assume that there exists a lift $\overline{F}$ of $F$ fitting into a commutative diagram as in Theorem \ref{t:podnoszenie} (2). The homomorphism $\mathrm{Cl}(Y) \rightarrow \mathrm{Cl}(X)$ that is a part of the data of the graded homomorphism of graded rings $\overline{F}^*:S\rightarrow R$ will be denoted by $\varphi$.

Let $Z$ be a MDS with the Cox ring $T$. The following simple example shows that non-isomorphic $\mathrm{Cl}(Z)$-graded $T$-modules can determine isomorphic quasicoherent sheaves on $Z$.
\begin{example}\label{p:rozne_moduly_ten_sam_snop}
Let $Z= \mathbb{P}_{\mathbb{C}}^1$. Then the Cox ring is $T=\mathbb{C}[x,y]$, and $\widehat{Z} = \mathbb{C}^2 \setminus \{(0,0)\}$. Let $M$ be the base field $\mathbb{C}$ with the structure of a $\mathbb{Z}$-graded $\mathbb{C}[x,y]$-module given by $x\alpha = y\alpha = 0$ for every $\alpha \in \mathbb{C}.$ Then $\overline{M}$ is a skyscraper sheaf on $\mathbb{C}^2$ supported at the origin. Hence ${i_Z}^*\overline{M} = 0$ and therefore $\widetilde{M} \cong \widetilde{0}$.
\end{example}
The above example suggests that we should make a choice of a particular $\mathrm{Cl}(Z)$-graded $T$-module describing a given quasicoherent sheaf $\mathcal{F}$ on $Z$. We will denote by $\Gamma_*(\mathcal{F})$ the $\mathrm{Cl}(Z)$-graded $T$-module $\Gamma(\widehat{Z}, {\pi_Z}^*\mathcal{F})$.

\subsection{The inverse image}\label{s:tii}
Let $N$ be a $\mathrm{Cl}(Y)$-graded $S$-module. Then $N \otimes_S R$ has a structure of an $R$-module. We will define its $\mathrm{Cl}(X)$-grading as follows: for homogeneous $n \in N$ and homogeneous $r \in R$ define ${\rm{deg}}(n\otimes r)=\varphi({\rm{deg}}(n))+{\rm{deg}}(r)$. It is straightforward to verify that this grading is well defined and gives $N \otimes_S R$ a structure of a $\mathrm{Cl}(X)$-graded $R$-module.

\begin{theorem}\label{t:pullback} In the setup from the beginning of the section, let $\mathcal{G}$ be a quasicoherent sheaf on $Y$. Assume that $\Gamma_*(\mathcal{G})\cong N$ for a $\mathrm{Cl}(Y)$-graded $S$-module $N$. Then $F^* \mathcal{G} \cong \widetilde{N\otimes_SR}$.
\end{theorem}
\begin{proof}
We are interested only in quasicoherent sheaves up to isomorphism, so it is enough to prove the theorem for $\mathcal{G} = \widetilde{N}$. From the commutativity of the diagram in Theorem \ref{t:podnoszenie} we have ${\pi_X}^*{F}^*\widetilde{N} = {\widehat{F}}^*{\pi_Y}^*{\widetilde{N}}$. Proposition $\ref{s:lematy_do_cofania}$ implies that ${\widehat{F}}^*{\pi_Y}^*{\widetilde{N}} \cong \widehat{F}^*{i_Y}^*{\overline{N}}$. Using once more the diagram in Theorem $\ref{t:podnoszenie}$ we obtain $\widehat{F}^*{i_Y}^*{\overline{N}} = {i_X}^*{\overline{F}}^*{\overline{N}}$. From the description of the inverse image of a quasicoherent sheaf by a morphism of affine schemes we obtain ${i_X}^* {\overline{F}}^*{\overline{N}} \cong {i_X}^*(\overline{N\otimes_SR})$.
Hence we have:
\[
(F^*\widetilde{N}) \otimes_{\mathcal{O}_X} \mathcal{R} \stackrel{\ref{s:lematy_do_cofania}}{\cong} {\pi_X}_*{\pi_X}^*{F^*{\widetilde{N}}} \cong {\pi_X}_*{i_X}^*(\overline{N \otimes_SR}).
\]
Taking the zeroth gradation we obtain an isomorphism $F^*{\widetilde{N}} \cong \widetilde{N\otimes_SR}.$
\end{proof}

\begin{remark}
In the notation from the above theorem, let $N'$ be a $\mathrm{Cl}(Y)$-graded $S$-module such that $\widetilde{N'}\cong \mathcal{G}$.
As shown in section 1.1.4 in \cite{BB13}, even for toric varieties, in general we do not have $\widetilde{N'\otimes_S R}\cong F^*\mathcal{G}$. However, we will show in Lemma \ref{l:dla_gladkich_moge_wziac_dowolny}, that if $Y$ is smooth, then these quasicoherent $\mathcal{O}_X$-modules are isomorphic. 
\end{remark}

\subsection{The direct image}
As we have seen in the beginning of Section \ref{s:tii}, the extension of scalars of a graded module by a graded homomorphism of graded rings gives a graded module. For the restriction of scalars it is not the case. Moreover, even if taking the restriction of scalars of a $\mathrm{Cl}(X)$-graded $R$-module $M$ gives a $\mathrm{Cl}(Y)$-graded $S$-module, it may be the case that it does not correspond to the direct image of $\widetilde{M}$ as the following example shows.
\begin{example}\label{p:istotnosc}
Let $Z=\mathbb{P}_{\mathbb{C}}^1$, $Z'=\mathbb{C}^1$ with coordinates $x,y$ and $t$, respectively. Let $F([x:y])=0.$ Consider the structure sheaf $\mathcal{O}_{\mathbb{P}_{\mathbb{C}}^1} \cong \widetilde{\mathbb{C}[x,y]}$. We have $\Gamma(Z', F_*\mathcal{O}_{\mathbb{P}_{\mathbb{C}}^1})=\mathbb{C}$ but $\Gamma(Z', \overline{F}_*\overline{\mathbb{C}[x,y]})=\mathbb{C}[x,y]$.
\end{example}
Let $M$ be a $\mathrm{Cl}(X)$-graded $R$-module. Let $M^*_S = \bigoplus_{[E] \in \mathrm{Cl}(Y)} M_{\varphi([E])}$, where $\varphi$ was defined at the beginning of Section \ref{s_3}. The graded homomorphism of graded rings $\overline{F}^*: S \rightarrow R$ gives $M^*_S$ a structure of a $\mathrm{Cl}(Y)$-graded $S$-module: for all $s\in S$ and for all $m \in M^*_S$ we define $s \cdot m=\overline{F}^*(s)\cdot m$.

\begin{theorem}\label{t:pushforward} In the setup from the beginning of the section, let $\mathcal{F}$ be a quasicoherent sheaf on $X$ with $\Gamma_*(\mathcal{F}) = M.$ Then $F_* \mathcal{F} \cong \widetilde{M^*_S}$.
\end{theorem}

In the proof we will define isomorphisms of sections of these two sheaves on a basis for the topology of $Y$. In order to be able to glue these isomorphisms to an isomorphism of quasicoherent sheaves we will carefully show that all isomorphisms considered on the way are natural. Before giving the proof of this theorem we will establish a few lemmas.

\begin{lemma}\label{l:Hartogs_type} Let $Z$ be a MDS with the Cox ring $T$. Let $\mathcal{F}$ be a quasicoherent sheaf on $Z$ with $\Gamma_*(\mathcal{F}) = M$ and let $g$, $h$ be homogeneous elements of $T$. Then there are commutative diagrams:
\begin{center}
\begin{tikzcd}[row sep=normal, column sep=tiny]
T_g \arrow{r} \arrow{d} & \mathcal{O}_{\overline{Z}}(\overline{Z}_g) \arrow{d} \arrow{r} &\mathcal{O}_{\widehat{Z}}(\widehat{Z}_g)=\mathcal{O}_{\overline{Z}}(\widehat{Z}_g) \arrow{d}\\
T_{hg} \arrow{r}  & \mathcal{O}_{\overline{Z}}(\overline{Z}_{hg}) \arrow{r} &\mathcal{O}_{\widehat{Z}}(\widehat{Z}_{hg}) = \mathcal{O}_{\overline{Z}}(\widehat{Z}_{hg})\\ 
T_h \arrow{r} \arrow{u} & \mathcal{O}_{\overline{Z}}(\overline{Z}_h) \arrow{u} \arrow{r} & \mathcal{O}_{\widehat{Z}}(\widehat{Z}_h)=\mathcal{O}_{\overline{Z}}(\widehat{Z}_h) \arrow{u}
\end{tikzcd}
\begin{tikzcd}[row sep=normal, column sep=tiny]
M_g \arrow{r} \arrow{d} & \Gamma(\overline{Z}_g, \overline{M})\arrow{d} \arrow{r} &\Gamma(\widehat{Z}_g, \overline{M}) \arrow{d}\\
M_{hg} \arrow{r}  & \Gamma(\overline{Z}_{hg}, \overline{M}) \arrow{r} &\Gamma(\widehat{Z}_{hg}, \overline{M})\\ 
M_h \arrow{r} \arrow{u} & \Gamma(\overline{Z}_{h}, \overline{M}) \arrow{u} \arrow{r} & \Gamma(\widehat{Z}_h, \overline{M}) \arrow{u}
\end{tikzcd}
\end{center}
with all horizontal arrows isomorphisms. In particular, for every homogeneous $h\in T$ there are isomorphisms $\alpha_h : T_h \rightarrow \mathcal{O}_{\widehat{Z}}(\widehat{Z}_h)$ and $\beta_h: M_h \rightarrow \Gamma(\widehat{Z}_h, \overline{M}).$
\end{lemma}

\begin{proof} Since $\overline{Z}$ is normal and the complement of $\widehat{Z}$ is of codimension at least two, it follows that restricting functions gives an isomorphism $\mathcal{O}_{\overline{Z}}(\overline{Z}) \xrightarrow{\cong} \mathcal{O}_{\widehat{Z}}(\widehat{Z})$. Hence we have $\mathcal{O}_{\overline{Z}} \cong {i_Z}_* \mathcal{O}_{\widehat{Z}}$. Therefore for every $g \in T$, restricting sections is an isomorphism $\mathcal{O}_{\overline{Z}}(\overline{Z}_g) \cong \mathcal{O}_{\widehat{Z}}(\widehat{Z}_g)$. This proves that the three right horizontal arrows of the left diagram are isomorphisms. It is well known that there exist three left horizontal arrows in this diagram, that are isomorphisms such that the left two squares commute (\cite{Har77} Proposition II.2.2). The right two squares commute since all maps are restrictions of the sections of the structure sheaf $\mathcal{O}_{\overline{Z}}$.

Since $\Gamma_*(\mathcal{F})=M$, $\overline{M} \cong {i_Z}_*{\pi_Z}^*\mathcal{F}$ as sheaves of abelian groups. Therefore for every non-zero homogeneous $f\in T$ the restriction of sections gives an isomorphism $\Gamma(\overline{Z}_f, \overline{M}) \rightarrow \Gamma(\widehat{Z}_f, \overline{M})$. Similar argument to the given above, shows that these isomorphisms give a commutative diagram as in the statement of the lemma.
\end{proof}
Given a surjective map of sets $G:Z_1 \rightarrow Z_2$, we say that a subset $U \subset Z_1$ is \textbf{saturated} with respect to $G$ if $G^{-1}G(U)=U.$

\begin{lemma}\label{l:natrualnosc_funkcji} Let $Z$ be a Mori Dream Space with the Cox ring $T$. Let $f,g$ be two (possibly zero) homogeneous elements of $T$ such that $\widehat{Z}_f$ and $\widehat{Z}_g$ are saturated with respect to $\pi_Z$. Then we have the following commutative diagram with obvious vertical maps:
\begin{center}
\begin{tikzcd}[row sep=normal, column sep=normal]
T_{(f)} \arrow{r}{\cong} \arrow{d} & \mathcal{O}_Z(\pi_Z(\widehat{Z}_f)) \arrow{d} \\
T_{(fg)} \arrow{r}{\cong}  & \mathcal{O}_Z(\pi_Z(\widehat{Z}_{fg}))\\
T_{(g)} \arrow{r}{\cong} \arrow{u} & \mathcal{O}_Z(\pi_Z(\widehat{Z}_g)) \arrow{u}
\end{tikzcd}
\end{center}
\end{lemma}

\begin{proof}
We have $({\pi_Z}_*\mathcal{O}_{\widehat{Z}})_0 \cong \mathcal{O}_Z$. Hence we have isomorphisms:
%
\[\Gamma(\pi_Z(\widehat{Z}_f),\mathcal{O}_{Z}) \cong \Gamma(\pi_Z(\widehat{Z}_f), {\pi_Z}_*\mathcal{O}_{\widehat{Z}})_0 = \Gamma({\pi_Z}^{-1}
(\pi_Z(\widehat{Z}_f)), \mathcal{O}_{\widehat{Z}})_0
 = \Gamma(\widehat{Z}_f, \mathcal{O}_{\widehat{Z}})_0 \cong T_{(f)}.
\]
%
Since the first isomorphism comes from the isomorphism of quasicoherent sheaves and the last comes from Lemma \ref{l:Hartogs_type}, they commute with restrictions and we have the commutative diagram from the statement.
We have used here a fact that intersection of saturated sets is saturated.
\end{proof}

\begin{proof}[Proof of Theorem \ref{t:pushforward}]
For all $[D]\in \mathrm{Cl}(Y)$ and for all non-zero $f \in {S}_{[D]}$ with $Y_{[D],f}$ affine we will define an isomorphism of $\mathcal{O}_Y(Y_{[D],f})$-modules: 
\[
\Gamma(Y_{[D],f},F_*\mathcal{F}) \xrightarrow{\chi_{[D],f}} \Gamma(Y_{[D],f}, \widetilde{M^*_S})
\]
such that for every $[E]\in \mathrm{Cl}(Y)$  and for every non-zero $g \in S_{[E]}$ we have the following commutative diagram:
\begin{center}
\begin{tikzcd}[row sep=normal, column sep=large]
\Gamma(Y_{[D],f},F_*\mathcal{F}) \arrow{d} \arrow{r}{\chi_{[D],f}} & \Gamma(Y_{[D],f}, \widetilde{M^*_S}) \arrow{d}\\
\Gamma(Y_{[D]+[E],fg},F_*\mathcal{F}) \arrow{r}{\chi_{[D]+[E],fg}} & \Gamma(Y_{[D]+[E],fg}, \widetilde{M^*_S})\\
\Gamma(Y_{[E],g},F_*\mathcal{F}) \arrow{u} \arrow{r}{\chi_{[E],g}} & \Gamma(Y_{[E],g}, \widetilde{M^*_S}) \arrow{u}
\end{tikzcd}
\end{center}
where the vertical arrows are restriction maps. By Lemmas \ref{l:baza_sklejenie_2} and \ref{l:baza_sklejanie}, it will follow that such maps $\chi_{[D],f}$ define an isomorphism of $\mathcal{O}_Y$-modules $F_*\mathcal{F} \rightarrow \widetilde{M^*_S}.$ Note that the identifications that we have already done in the lemmas are all natural in the sense that they fit into similar diagrams.

\textbf{Step 1.} Pick any $[D]\in \mathrm{Cl}(Y)$ and any non-zero $f \in S_{[D]}$ such that $Y_{[D],f}$ is affine. By Proposition 1.6.3.3 in \cite{ADHL15}, we have $\pi_Y^{-1}(Y_{[D],f}) = \widehat{Y}_f=\overline{Y}_f$. Therefore, since $\pi_Y$ is surjective, $\widehat{Y}_f$ is saturated with respect to $\pi_Y$. Hence by Lemma \ref{l:natrualnosc_funkcji} we may assume that for every $[D]\in \mathrm{Cl}(Y)$ and for every non-zero $f \in S_{[D]}$ such that $Y_{[D],f}$ is affine we have $\mathcal{O}_{Y}(Y_{[D],f}) = S_{(f)}.$
We will describe $F_*\mathcal{F}|_{Y_{[D],f}}$. Since $Y_{[D],f}$ is affine it is enough to compute $\Gamma(Y_{[D],f},F_*\mathcal{F})$ and describe its $S_{(f)}$-module structure.

Since $Y_{[D],f}$ is affine, $\pi_Y^{-1}(Y_{[D],f}) = \widehat{Y}_f = \overline{Y}_f.$ We have also $\widehat{F}^{-1}(\widehat{Y}_f) = \widehat{X}_{f\circ \overline{F}}.$ From this equality we obtain from the diagram in Theorem \ref{t:podnoszenie} that: 
\[
\pi_X^{-1}(F^{-1}(Y_{[D],f})) = \widehat{X}_{f\circ \overline{F}}.
\]
Hence from surjectivity of $\pi_X$ it follows that $F^{-1}(Y_{[D],f}) = \pi_X(\widehat{X}_{f\circ \overline{F}})$ and we have the following commutative diagram:
\begin{center}
\begin{tikzcd}[row sep=normal, column sep=normal]
\overline{X}_{f\circ \overline{F}} \arrow[dotted]{r}{\overline{F}} & \overline{Y}_f \\
\widehat{X}_{f\circ \overline{F}} \arrow[hook]{u}{i_X} \arrow[dotted]{r}{\widehat{F}} \arrow[two heads]{d}{\pi_X} &
\widehat{Y}_f \arrow[two heads]{d}{\pi_Y} \arrow{u}{=} \\
\pi_X(\widehat{X}_{f\circ \overline{F}}) \arrow{r}{F} 	& Y_{[D],f}
\end{tikzcd}
\end{center}

It follows that for every  $[D]\in \mathrm{Cl}(Y)$ and for every non-zero $f \in S_{[D]}$ such that $Y_{[D],f}$ is affine, $\widehat{X}_{f\circ \overline{F}}$ is saturated with respect to $\pi_X$ and therefore by Lemma \ref{l:natrualnosc_funkcji} we may assume that for such $[D]$ and $f$ we have $\mathcal{O}_X(\pi_X(\widehat{X}_{f\circ \overline{F}}))=R_{(f\circ \overline{F})}.$

By Proposition \ref{s:lematy_do_cofania} $\mathcal{F}\otimes_{\mathcal{O}_X}\mathcal{R} = {\pi_X}_*{\pi_X}^*\mathcal{F}$ so  $\Gamma(\pi_X(\widehat{X}_{f\circ \overline{F}}), \mathcal{F})$ is the degree zero part of $\Gamma(\widehat{X}_{f \circ \overline{F}}, {\pi_X}^*\mathcal{F})$ which by naturality of ${\beta_h}'s$ in Lemma \ref{l:Hartogs_type} can be assumed to be equal to $M_{f\circ \overline{F}}$. We have established that $\Gamma({Y_{[D],f}}, F_*\mathcal{F}) = M_{(f\circ \overline{F})}$ therefore describing the group structure of $\Gamma(Y_{[D],f},F_*\mathcal{F})$.

\textbf{Step 2.} We will now describe the $S_{(f)}$-module structure of $\Gamma(Y_{[D],f},F_*\mathcal{F})$.
Firstly we want to describe the module structure on $\Gamma(F^{-1}(Y_{[D],f}),\mathcal{F})$. 
From the diagram in Theorem \ref{t:podnoszenie} and the description of quasicoherent sheaves on Mori Dream Spaces from Proposition \ref{s:funktor} we know that $\Gamma(F^{-1}(Y_{[D],f}),\mathcal{F})$ is the degree zero part of $\Gamma(F^{-1}(Y_{[D],f}),{\pi_X}_*{i_X}^*\overline{M}) = \Gamma(\widehat{X}_{f\circ \overline{F}},{i_X}^*\overline{M})$. Hence it is $M_{(f\circ \overline{F})}$ with the $\mathcal{O}_X(\pi_X(\widehat{X}_{f \circ \overline{F}}))=R_{(f\circ \overline{F})}$-module structure coming from the map $\mathcal{O}_X(\pi_X(\widehat{X}_{f \circ \overline{F}})) \rightarrow {\pi_X}_*\mathcal{O}_{\widehat{X}}(\pi_X(\widehat{X}_{f \circ \overline{F}})) = \mathcal{O}_{\widehat{X}}(\widehat{X}_{f\circ \overline{F}}).$ This map is the inclusion $R_{(f\circ \overline{F})} \rightarrow R_{f\circ\overline{F}}$. Hence $\Gamma(F^{-1}(Y_{[D], f}), \mathcal{F})$ is $M_{(f\circ \overline{F})}$ not only as an abelian group but also as an $R_{(f \circ \overline{F})}$-module. Therefore $\Gamma(Y_{[D], f}, F_*\mathcal{F}) = M_{(f \circ \overline{F})}$ with the $S_{(f)}$-module structure coming from the map $\mathcal{O}_Y(Y_{[D], f}) \rightarrow F_*\mathcal{O}_X(Y_{[D],f}).$ Which is the map $S_{(f)} \rightarrow R_{(f\circ \overline{F})}$. Therefore, up to natural isomorphisms, $\Gamma(Y_{[D], f}, F_*\mathcal{F})=M_{(f\circ \overline{F})}$ as an $S_{(f)}$-module.

\textbf{Step 3.} We will describe the sections of $\widetilde{M^*_S}$ over affine sets of the form $Y_{[D],f}$. We will assume, using the naturality of ${\alpha_h}'s$ in Lemma \ref{l:Hartogs_type}, that $\mathcal{O}_{\widehat{Y}}(\widehat{Y}_f) = S_f$ and using Lemma \ref{l:natrualnosc_funkcji} that $\mathcal{O}_{Y}(Y_{[D],f})=S_{(f)}$. Then from the description of quasicoherent sheaves on Mori Dream Spaces we have $\Gamma(Y_{[D],f}, \widetilde{M^*_S})=({M^*_S})_{(f)}$ as an $S_{(f)}$-module. We have ${\rm{deg}}(f)=[D]$ and ${\rm{deg}}(f\circ \overline{F})=\varphi([D])$ hence: 
\[
M_{(f\circ \overline{F})} = \{\frac{m}{(f\circ \overline{F})^n} | n\in \mathbb{N},\  m\in M \textrm{ and } {\rm{deg}}(m) = n\varphi([D])\}
\]
and 
\[
(M^*_S)_{(f)} = \{\frac{m}{{f}^n} | n \in \mathbb{N},\ m \in M^*_S \textrm{ and } {\rm{deg}}(m) = n[D]\}.
\]
From the definition of $M^*_S$ it follows that we have $(M^*_S)_{k[D]} = M_{k\varphi([D])}$ as abelian groups. Therefore: 
\[
(M^*_S)_{(f)} = \{\frac{m}{{f}^n} | n \in \mathbb{N},\ m \in M \textrm{ and }{\rm{deg}}(m) = n\varphi([D])\}
\]
and we have an isomorphism of $S_{(f)}$-modules $\chi_{[D],f}: M_{(f \circ \overline{F})} \rightarrow (M^*_S)_{(f)}$ given by $\frac{m}{(f \circ \overline{F})^k} \mapsto \frac{m}{f^k}$ for $m \in M_{k\varphi([D])}.$ This isomorphism is natural so isomorphisms of this type for all affine sets of the form $Y_{[D],f}$ will glue by Lemma \ref{l:baza_sklejenie_2} to an isomorphism of $\mathcal{O}_Y$-modules $F_* \mathcal{F} \rightarrow \widetilde{M^*_S}.$  Observe that $\chi_{[D],f}$ is well defined. If $\frac{m}{(f \circ \overline{F})^k}=\frac{n}{(f \circ \overline{F})^l}$, then there exists $s\in \mathbb{N}$, such that $(f \circ \overline{F})^s((f \circ \overline{F})^lm-(f \circ \overline{F})^kn)=0$. Then by definition of the $S_{(f)}$-module structure on $(M^*_S)_{(f)}$ we have $f^s(f^lm-f^kn)=0$.

\end{proof}

\begin{remark} In the notation from Examples \ref{p:rozne_moduly_ten_sam_snop} and \ref{p:istotnosc} we have $\widetilde{M}\cong 0$ but $M_S^* = M_0 =\mathbb{C}$. 
Thus $\widetilde{M^*_S}$ is not isomorphic to $F_*\widetilde{M}$. Therefore, in general, in Theorem \ref{t:pushforward} we cannot use arbitrary $\mathrm{Cl}(X)$-graded $R$-module $M$ such that $\widetilde{M} \cong \mathcal{F}$.
\end{remark}

\subsection{Additional remarks and special cases}
We start with describing the module $\Gamma_*(\mathcal{O}_X(D))$ for a line bundle $\mathcal{O}_X(D)$.

\begin{lemma}\label{l:gammadladywizorow} Let $X$ be a smooth MDS with the Cox ring $R$. Given a Cartier divisor $D\in \mathrm{Pic}(X)$ we have $\Gamma_*(\mathcal{O}_X(D))\cong R([D])$.
\end{lemma}
\begin{proof} In the notation from Construction \ref{k:iloraz} we have $\Gamma_*(\mathcal{O}_X(D))=\Gamma(\widehat{X}, {\pi_X}^*(\mathcal{O}_X(D)))$.
Since $\pi_X$ is surjective and $D$ is Cartier we have ${\pi_X}^*(\mathcal{O}_X(D))=\mathcal{O}_{\widehat{X}}(\pi^*D)$ where $\pi^*D$ is the pullback of the divisor $D$. Write $D$ as $D=E_1 - E_2$ where both $E_1$ and $E_2$ are effective. From Proposition 1.5.2.2 in \cite{ADHL15} there exist $[D_1]$, $[D_2]\in \mathrm{Cl}(X)$ and $f_1\in R_{[D_1]}$, $f_2\in R_{[D_2]}$ such that $E_i = {\rm{div}}_{[D_i]}(f_i)$. From Lemma 1.5.3.6 in \cite{ADHL15} it follows that ${\pi_X}^*D = {\rm{div}}(\frac{f_1}{f_2})$. Hence $\Gamma_*(\mathcal{O}_X(D))\cong R([D_1]-[D_2])$. By the definition of the $[D]$-divisor, $[D_1]-[D_2]=[D].$
\end{proof}

If we add additional assumptions to the MDSes that we consider and the quasicoherent sheaves on them, the situation is easier.
It is shown by the following lemmas that are all easy consequences of Proposition 1.6.1.6 in \cite{ADHL15} and  Theorem 4.2.14 in \cite{HL10}.

\begin{lemma}\label{l:dla_gladkich_moge_wziac_dowolny}
In the setting of Theorem \ref{t:pullback}, let $Y$ be smooth and let $\mathcal{G}$ be coherent. If $N'$ is any $\mathrm{Cl}(Y)$-graded $S$-module such that $\widetilde{N'}\cong \mathcal{G}$, then $\widetilde{N'\otimes_S R} \cong F^*\mathcal{G}$.
\end{lemma}
\begin{proof}
If $Y$ is smooth, $\pi_Y:\widehat{Y}\to Y$ is a principal $H_Y=\operatorname{Spec}k[\mathrm{Cl}(X)]$-bundle by Proposition 1.6.1.6 in \cite{ADHL15}. Hence, by Theorem 4.2.14 in \cite{HL10} and Proposition \ref{s:lematy_do_cofania},  $\widetilde{N'} \cong \mathcal{G}$ implies $\pi_Y^*\mathcal{G} \cong i_Y^*(\overline{N'})$. Therefore we have an isomorphism $i_X^*\overline{N'\otimes_S R} \cong \pi_X^* F^*\mathcal{G}$. The claim follows from Propositon \ref{s:lematy_do_cofania}.
\end{proof}

To simplify notation, in the next two lemmas we will write $f\sim g$ if $f:A\to B$, $g:C\to D$ are morphisms in a category and there exist isomorphisms $a:A\to C$, $b:B\to D$ such that $g\circ a = b\circ f$.

\begin{lemma}\label{l:gladkieMDS_morfizmy_pochodza_z_algebry}
Let $X$ be a smooth MDS with the Cox ring $R$. Let $\varphi:\mathcal{F}\to\mathcal{G}$ be a homomorphism of coherent sheaves on $X$. Then there exists a graded homomorphism of $\mathrm{Cl}(X)$-graded $R$-modules $\alpha : M\to N$ such that $\widetilde{\alpha} \sim \varphi$. Moreover, we can take $\alpha$ to be the homomorphism obtained by applying $\Gamma(\widehat{Y},.)$ to $\pi_Y^*(\varphi)$.
\end{lemma}
\begin{proof}
Let $\alpha$ be the graded homomoprhism of $\mathrm{Cl}(X)$-graded $R$-modules obtained by 
applying $\Gamma(\overline{Y},.)$ to the morphism of quasicoherent sheaves $(i_Y)_* \pi_Y^* 
\varphi$. We have the following relations: 
\begin{itemize}
\item[i)]$\overline{\alpha} \sim (i_Y)_* \pi_Y^*\varphi$,
\item[ii)]$i_Y^* (i_Y)_*(\pi_Y^*\varphi)  \sim (\pi_Y^*\varphi)$,
\item[iii)] $((\pi_Y)_*\pi_Y^* \varphi)_0 \sim \varphi$.
\end{itemize}
The first two relations are obvious. The third follows from Theorem 4.2.14 in \cite{HL10}. All three imply that $\widetilde{\alpha}\sim \varphi$.
\end{proof}

\begin{lemma}\label{l:morfizm_dla_wiazek}
Let $X$ be a smooth MDS with the  Cox ring $R$. If $\alpha, \beta:\oplus_{i=1}^n R(a_i)\to \oplus_{j=i}^m R(b_i)$ are graded homomorphisms of $\mathrm{Cl}(X)$-graded $R$-modules such that $\tilde{\alpha}=\tilde{\beta}$, then $\alpha \sim \beta$.
\end{lemma}
\begin{proof}
If $\tilde{\alpha} = \tilde{\beta}$, then by Theorem 4.2.14 in \cite{HL10} we have $(i_X)^*\overline{\alpha} \sim (i_X)^*\overline{\beta}$. Since $(i_X)_* \mathcal{O}_{\widehat{X}} \cong \mathcal{O}_{\overline{X}}$ and since we are dealing with locally free $\mathcal{O}_{\overline{X}}$-modules of finite rank, we can use the projection formula to obtain:
\[
\overline{\alpha} \sim (i_X)_*i_X^*\overline{\alpha} \sim (i_X)_*i_X^*\overline{\beta}\sim \overline{\beta}.
\]
\end{proof}

\section{Examples}
In this section we will present two examples. We will work with varieties over $\mathbb{C}$ since this is the assumption made in the book \cite{CLS11}, from which we will cite some results. 

The first example will be a pushforward of tangent sheaf under a toric morphism of smooth toric varieties. This example will show, that in the notation from Theorem \ref{t:pushforward}, in general $\Gamma_*(F_*\mathcal{F})$ is not isomorphic to $(\Gamma_*(\mathcal{F}))_S^*$. The second example is slightly more complicated since the target MDS is not a toric variety. 

Both examples illustrate that the main problem in applying Theorems \ref{t:pullback} and \ref{t:pushforward} is finding the right graded module describing a given quasicoherent sheaf. We skip most of the verifications. 
Some calculations were done in Macaulay2 \cite{M2}. 

\subsection{Tangent sheaf of the Hirzebruch surface}
We will consider the pushforward of the tangent sheaf of the Hirzebruch surface under the toric morphism to $\mathbb{P}_\mathbb{C}^1$ induced by the projection $\mathbb{R}^2 \rightarrow \mathbb{R}$ onto the $x$-axis.

\begin{example} Let $M,N$ be dual lattices of rank two. Fix a natural number $a\in \mathbb{N}$. Let $X$ be the Hirzebruch surface $\mathbb{F}_a$ given by the unique complete fan $\Sigma_1$ in $N_\mathbb{R} := N\otimes_\mathbb{Z}\mathbb{R} \cong \mathbb{R}^2$ with ray generators $u_{\rho_1} = (1,0)$, $u_{\rho_2}=(0,1)$, $u_{\rho_3}=(-1,-a)$ and $u_{\rho_4}=(0,-1)$. Let $Y$ be the projective line given 
by the unique complete fan $\Sigma_2$ in $\mathbb{R}$. Denote the ray generators of $\Sigma_2$ by $w_1=1$ and $w_2=-1$. Let $\phi: N \cong \mathbb{Z}^2 \rightarrow \mathbb{Z}$ be given by $(x_1,x_2)\mapsto x_1$. The tensored map $\phi_\mathbb{R}:N_\mathbb{R} \rightarrow \mathbb{R}$ is compatible with the fans $\Sigma_1$ and $\Sigma_2$. Hence it induces a toric morphism $F: X \rightarrow Y$. We will denote the Cox rings of $X$ and $Y$ by $R$ and $S$, respectively. In the following description of the Cox rings of $X$ and $Y$ we use results from section  \S 5.1 in \cite{CLS11}.

The class group of $X$ is isomorphic to $\mathbb{Z}^2$. The Cox ring is given by $R=\mathbb{C}[x_1,x_2, x_3, x_4]$ with ${\rm{deg}}(x_1)={\rm{deg}}(x_3)=(1,0)$, ${\rm{deg}}(x_2)=(a,1)$ and ${\rm{deg}}(x_4)=(0,1)$. The Cox ring of $Y$ is $S=\mathbb{C}[y_1,y_2]$ with ${\rm{deg}}(y_1)={\rm{deg}}(y_2)=1.$
The irrelevant ideal of $X$ is generated by $x_1x_2$, $x_2x_3$, $x_3x_4$ and $x_1x_4$. The torus action on $X$ is given by $(\lambda, \mu)(x_1,x_2,x_3,x_4)=(\lambda x_1, \lambda^a\mu x_2, \lambda x_3, \mu x_4)$.

Let $\overline{F}:\mathbb{C}^4\rightarrow \mathbb{C}^2$ be given by $(s,t,u,w)\mapsto (s,u)$. On the level of coordinate rings it is given by $\overline{F}^*:\mathbb{C}[y_1,y_2]\rightarrow \mathbb{C}[x_1,x_2,x_3,x_4]$, where $y_1\mapsto x_1$ and $y_2 \mapsto x_3$. Thus it is a graded homomorphism of graded rings. We will show that it restricts to a map $\widehat{F}:\widehat{X}\rightarrow \widehat{Y}$. Suppose that $\overline{F}(s,t,u,w)=(0,0)$. Then $s=u=0$. Hence $(s,t,u,w)\in \overline{X}\setminus \widehat{X}$. It can be checked on the affine open covers associated with maximal cones, that the morphism $X\to Y$ induced by $\overline{F}$ is equal to $F$.

Let $\mathcal{T}_X$ be the tangent sheaf of $X$. Let $D_{\rho_i}$ be the divisor associated with the cone $\operatorname{cone}(\rho_i)$ for $i=1,...,4$.  Consider the map $\alpha: R\oplus R \rightarrow \bigoplus_{i=1}^4 R([D_{\rho_i}])$ given by $(s,t)\mapsto (x_1s, x_2(as+t),x_3s,x_4t)$. It is clearly an injective homomorphism of $\mathrm{Cl}(X)$-graded $R$-modules. Let $P$ denote the cokernel of this map. We have an exact sequence:
\begin{equation}\label{e:ciagdlastycznej}
0 \rightarrow R\oplus R \xrightarrow{\alpha} \bigoplus_{i=1}^4 R([D_{\rho_i}]) \to P \rightarrow 0.
\end{equation} 

We will later show in Lemma \ref{l:wiazkastyczna} that $\widetilde{P^*_S} \cong ((\Gamma_*(\mathcal{T}_X))^*_S)$ $\widetilde{ }$ $ $, where $ $ $\widetilde{ }$ $ $ is the functor from Proposition \ref{s:funktor}, $\Gamma_*$ and $(-)^*_S$ were defined in Section \ref{s_3}. That is $P$ can be used to compute the pushforward of $\mathcal{T}_X$ using Theorem \ref{t:pushforward}. Assuming this fact, we will describe the direct image sheaf.

The map $\phi: \mathbb{Z} \cong \mathrm{Cl}(Y) \rightarrow \mathrm{Cl}(X) \cong \mathbb{Z}^2$ associated with $\overline{F}^* : S \rightarrow R$ is given by $n\mapsto (n,0)$. Thus, by Theorem \ref{t:pushforward} and Lemma \ref{l:wiazkastyczna}, $F_*\mathcal{T}_X$ is the quasicoherent sheaf associated with the $\mathbb{Z}$-graded $S$-module $\bigoplus_{n\in \mathbb{Z}} P_{n,0}$. It can be checked that the $\mathbb{Z}$-graded $S$-module $P^*_S = \bigoplus_{n\in \mathbb{Z}} P_{n,0}$ is isomorphic to $S(-a)\oplus S(1) \oplus S(1) \oplus S(a)$. Therefore, Theorem \ref{t:pushforward} and Lemma \ref{l:gammadladywizorow} imply that $F_* \mathcal{T}_X \cong \mathcal{O}(-a)\oplus \mathcal{O}(1)\oplus\mathcal{O}(1)\oplus\mathcal{O}(a)$.

We are left with the proof that $P$ can be used instead of $\Gamma_*(\mathcal{T}_X)$ to calculate $F_* \mathcal{T}_X$. 
For an $R$-module $M$, let $M^\vee$ be the $R$-module $\operatorname{Hom}_R(M,R)$. 
If $M$ is a $\mathrm{Cl}(X)$-graded $R$-module we have a submodule $\operatorname{Hom}^\bullet_R(M,R) = \bigoplus_{[D]\in \mathrm{Cl}(X)} \operatorname{Hom}_R^{[D]}(M,R) \subseteq M^\vee$, where $\operatorname{Hom}_R^{[D]}(M,R)$ are graded homomorphisms of degree $[D]$, i.e. $\delta(M_{[E]})\subseteq R_{[D]+[E]}$ for a morphism $\delta \in \operatorname{Hom}_R^{[D]}(M,R)$. If $M$ is a finitely generated $\mathrm{Cl}(X)$-graded $R$-module, then $\operatorname{Hom}^\bullet_R(M,R) = M^\vee$. Therefore $M^\vee$ is a $\mathrm{Cl}(X)$-graded $R$-module.

Let $\theta: \bigoplus_{i=1}^4 R(-[D_{\rho_i}]) \rightarrow R\oplus R$ be given by $(f_1,...,f_4)\mapsto (x_1f_1+ax_2f_2+x_3f_3, x_2f_2+x_4f_4)$. Let $Q$ be the kernel of $\theta$. It can be checked that we have $\Gamma_*(\Omega_X^1)\cong Q$ and $\Gamma_*(\mathcal{T}_X)\cong Q^\vee$. Moreover, applying the functor  
$\operatorname{Hom}_R(\cdot, R)$ to exact sequence (\ref{e:ciagdlastycznej}). We obtain the exact sequence:
\[
0 \rightarrow P^\vee \rightarrow \bigoplus_{i=1}^4 R(-[D_{\rho_i}]) \xrightarrow{\theta} R\oplus R.
\]
Hence $P^\vee \cong Q$ as both are the kernel of $\theta$.

\begin{lemma}\label{l:wiazkastyczna}
We have $\widetilde{P^*_S} \cong ((\Gamma_*(\mathcal{T}_X))^*_S)$ $\widetilde{ }$ $ $.
\end{lemma}
\begin{proof}
Let $h: P \rightarrow P^{\vee \vee}$ be the natural map. We claim that it is injective. Indeed, $P$ is torsion-free and it is finitely generated over $R$. Thus it is isomorphic to a submodule of a finitely generated free $R$-module. Hence the map $h$ is injective by Exercise 1.4.20 in \cite{BH98}. Moreover it is clearly a graded homomorphism of degree 
zero. Since $\Gamma_*(\mathcal{T}_X)\cong Q^\vee$, it is enough to show that for large enough $n$ we have equality of dimensions of $P^{\vee \vee}_{n,0} \cong Q^\vee_{n,0}$ and $P_{n,0}$ as $\mathbb{C}$-vector spaces.
We omit this easy calculation.
\end{proof}
\end{example}

\begin{remark} We showed that $F_* \mathcal{T}_X \cong \mathcal{O}(-a)\oplus \mathcal{O}(1)\oplus\mathcal{O}(1)\oplus\mathcal{O}(a)$. Therefore, by Lemma \ref{l:gammadladywizorow}, $\Gamma_*(F_*\mathcal{T}_X) \cong S(-a)\oplus S(1) \oplus S(1)\oplus S(a)$. We have $\Gamma_*(\mathcal{T}_X)=Q^\vee$. Therefore Theorem \ref{t:pushforward} gives the $\mathbb{Z}$-graded $S$-module $(Q^\vee)^*_S$ as the one describing $F_* \mathcal{T}_X$. 
Let $a=0$. Then,  $Q^\vee \cong R(2,0)\oplus R(0,2)$. Thus the graded part of $(Q^\vee)^*_S$ in degree $-2$ has dimension $1$ as a $\mathbb{C}$-vector space. However, $S \oplus S(1)\oplus S(1) \oplus S$ clearly has no non-zero homogeneous element of degree $-2$. Therefore, in the setting of Theorem \ref{t:pushforward}, we do not in general have isomorphisms of $\mathrm{Cl}(Y)$-graded $S$-modules $\Gamma_*(F_*\mathcal{F})$ and $(\Gamma_*(\mathcal{F}))^*_S$.

\end{remark}

\subsection{Smooth quintic del Pezzo surface}
Let $X$ be a smooth quintic del Pezzo surface, i.e. the blow-up of four points $p_1$, $p_2$, $p_3$ and $p_4$ in $\mathbb{P}^2$ such that no three of them lie on a line. The following description of the Cox ring of $X$ comes from section 5.2 in \cite{ADHL15}.

The class group of $X$ is given by $\mathrm{Cl}(X) =\mathbb{Z}[H]\oplus\mathbb{Z}[E_1]\oplus \mathbb{Z}[E_2]\oplus \mathbb{Z}[E_3]\oplus\mathbb{Z}[E_4]$, where $H$ is the strict transform of a general line in $\mathbb{P}^2$ and $E_i$ is the exceptional divisor over $p_i$. We identify $\mathrm{Cl}(X)$ with $\mathbb{Z}^5$ in the natural way, i.e. $[H]=e_1,[E_i]=e_{i+1}$ for $i=1,...,4$. For $1\leq i < j \leq 4$ let $L_{ij}$ denote the strict transform of the line in $\mathbb{P}^2$ passing through $p_i$ and $p_j$.
The Cox ring of $X$ is generated by canonical sections $x_i = 1\in \Gamma(X,\mathcal{O}_{X}(E_i))$ for $i=1,..., 4$ and $x_{ij} = 1 \in \Gamma(X,\mathcal{O}_X(L_{ij})$ for $1\leq i < j \leq 4$.
Let $S=\mathbb{C}[x_1,x_2,x_3,x_4,x_{12},x_{13},x_{14},x_{23},x_{24},x_{34}]$ be the $\mathbb{Z}^5$-graded polynomial ring with $deg(x_i)=e_{i+1}$ for $i=1,...,4$ and $deg(x_{ij})=e_1-e_{i+1}-e_{j+1}$ for $1\leq i < j \leq 4$.
The Cox ring of $X$ is the quotient of $S$ by the ideal $I=(x_2x_{12}-x_3x_{13}+x_4x_{14}, x_1x_{12}-x_3x_{23}+x_4x_{24}, x_1x_{13}-x_2x_{23}+x_4x_{34},x_1x_{14}-x_2x_{24}+x_3x_{34},x_{12}x_{34}-x_{13}x_{24}+x_{14}x_{23})$. We will denote the quotient $S/I$ by $R$.

Let $J = I + (x_{12}x_{34}+x_{13}x_{24})$. This ideal is prime in $S$ so it describes a one dimensional closed, irreducible subvariety $Y$ of $X$. 
We claim that $Y$ is the strict transform of a conic in $\mathbb{P}^2$ passing through $p_1,...,p_4$. To see this, note that $\operatorname{dim}_\mathbb{C}R_{(2,-1,-1,-1,-1)}=2$ is equal to the dimension of the vector space of homogeneous polynomials of degree $2$ that has zeroes at $p_1,...,p_4$. Therefore we have a pencil of conics in $\mathbb{P}^2$ passing through $p_1,...,p_4$ and their strict transforms are the only curves in $X$ defined by an element of $R$ of degree $(2,-1,-1,-1,-1)$.
We will calculate the restriction to $Y$ of the cotangent sheaf $\Omega_X^1$.

Let $A=\mathbb{C}[a,b]$ be the Cox ring of $\mathbb{P}^1$. Let $\overline{F}^*: S \to A$ be given by $x_1 \mapsto a+b$, $x_2 \mapsto 2a$, $x_3 \mapsto 2b$, $x_4\mapsto 2(a-b), x_{12}\mapsto 2, x_{13}\mapsto 2, x_{14} \mapsto -2, x_{23} \mapsto 2, x_{24}\mapsto -1, x_{34}\mapsto 1$.

Then $\overline{F}^*$ is a graded homomorphism of graded rings. It corresponds to an equivariant morphism $\overline{F}:\mathbb{A}^2\to \mathbb{A}^{10}$.
Using the anticanonical class of $X$ in Corrolary 1.6.3.6 in \cite{ADHL15} we can determine the irrelevant ideal of $X$. Then, it can be checked that $\overline{F}$ induces a morphism $F:\mathbb{P}^1\to X$ fitting into a commutative diagram analogous to the one in Theorem \ref{t:podnoszenie} (2). Moreover, $F:\mathbb{P}^1\to Y$ is an isomorphism.

We need to find a $\mathbb{Z}^5$-graded $R$-module $M$ such that $\widetilde{M} \cong \Omega_X^1$. We will embed $X$ into a product of projective spaces in order to be able to use the Euler sequence. Let $[x:y:z]$ denote the coordinates on $\mathbb{P}^2$ and let $[\alpha:\beta]$ denote the coordinates on $\mathbb{P}^1$. We will assume that the points that are blown-up are $p_1=[1:0:0]$, $p_2=[0:1:0]$, $p_3=[0:0:1]$ and $p_4=[1:1:1]$. 
Consider the rational map $\mathbb{P}^2 \dashrightarrow \mathbb{P}^1$ given by the pencil of conics passing through the points $p_1,...,p_4$. It is given by
$[x:y:z]\mapsto [xy-yz:yz-xz]$. The composition with the blowing-up $\pi: X \to \mathbb{P}^2$ gives a morphism $X\to \mathbb{P}^1$. This gives an embedding of $X$ into $\mathbb{P}^1\times \mathbb{P}^2$ that is onto $Z=V(\alpha(xz-yz)+\beta(xy-yz))$. 
Let $B=\mathbb{C}[\alpha,\beta,x,y,z]$ be a $\mathbb{Z}^2$-graded polynomial ring with $deg(\alpha)=deg(\beta) = e_1$ and $deg(x)=deg(y)=deg(z)=e_2$. Then $B$ is the Cox ring of $\mathbb{P}^1\times\mathbb{P}^2$. Consider the ring homomorphism $\overline{G}^*:B\to R$ given by:
$\alpha \mapsto \overline{x_{14}}\overline{x_{23}}$, $\beta \mapsto -\overline{x_{13}}\overline{x_{24}}$
$x\mapsto \overline{x_1}\overline{x_2}\overline{x_{12}}$, $y\mapsto \overline{x_2}\overline{x_3}\overline{x_{23}}$, $z\mapsto \overline{x_1}\overline{x_3}\overline{x_{13}}$. The homomorphism $\overline{G}^*$ is a graded homomorphism of graded rings. 
It corresponds to an equivariant map $\overline{G}:\overline{X} \to \mathbb{A}^5$. It can be checked that $\overline{G}$ induces a morphism $G:X\to \mathbb{P}^1\times \mathbb{P}^2$ fitting into a commutative diagram similar to the one in Theorem \ref{t:podnoszenie} (2). Moreover, $G:X\to Z$ is an isomorphism.

We will use the Euler sequence for the cotangent sheaf of projective space.
Let $C=\mathbb{C}[\alpha,\beta]$ be the Cox ring of $\mathbb{P}^1$.
Consider the graded homomorphism of $\mathbb{Z}$-graded $C$-modules $\gamma:C(-1)^{\oplus 2} \to C$ given by $(f,g) \mapsto \alpha f+ \beta g$. Let $M'$ denote the kernel of this map. Then $\widetilde{M'} \cong \Omega^1_{\mathbb{P}^1}$.
Similarly, let $D=\mathbb{C}[x,y,z]$ be the Cox ring of $\mathbb{P}^2$. Let $N'$ be the kernel of the graded homomorphism of $\mathbb{Z}$-graded $D$-modules $\delta: D^{\oplus 3}(-1)\to D$, given by $(f,g,h)\mapsto xf+yg+zh$. Then $\widetilde{N'} \cong \Omega^1_{\mathbb{P}^2}$.
We will denote by $M$ the $B$-module obtained from $M'$ by the extension of scalars along the inclusion $C\subset B$. We will do similarly for $N$.
Then $\widetilde{M\oplus N} \cong \Omega_{\mathbb{P}^1\times \mathbb{P}^2}^1$.
Moreover, since $\mathbb{P}^1\times\mathbb{P}^2$ is smooth, the morphism $\pi_{\mathbb{P}^1\times\mathbb{P}^2}$ is flat by Proposition 1.6.1.6 in \cite{ADHL15}. Therefore, the functor $\pi_{\mathbb{P}^1\times\mathbb{P}^2}^*$ is exact. It follows from Lemmas \ref{l:gammadladywizorow}, \ref{l:gladkieMDS_morfizmy_pochodza_z_algebry} and \ref{l:morfizm_dla_wiazek} that $\Gamma_*(\Omega_{\mathbb{P}^1\times \mathbb{P}^2}^1)$ is the kernel of the map $\epsilon = \gamma\otimes_C 1_B \oplus \delta\otimes_D 1_B$. In particular, it is isomorphic to $M\oplus N$.

Let $\mathcal{I}$ denote the ideal sheaf of $Z\subset \mathbb{P}^1\times \mathbb{P}^2$. Let $i:Z\to \mathbb{P}^1\times \mathbb{P}^2$ denote the inclusion.  Since $Z\subseteq \mathbb{P}^1\times \mathbb{P}^2$ is given by a single equation of bidegree (1,2), we have $\mathcal{O}_{\mathbb{P}^1\times \mathbb{P}^2}(-1,-2) \cong \mathcal{I}/\mathcal{I}^2$. We will find a $\mathrm{Cl}(X)$-graded $R$-module $Q$ that describes the cotangent sheaf $\Omega_X^1$ using the following diagram of coherent sheaves on $X$ with exact row and column.

\begin{center}
\begin{tikzcd}[row sep=small, column sep=small]
0\arrow{r} & \arrow{r} \mathcal{O}_X(-4,1,1,1,1) & G^*\Omega_{\mathbb{P}^1\times \mathbb{P}^2}^1 \arrow{r} \arrow{d} & \Omega_X^1 \arrow{r} & 0\\
 & & \mathcal{O}_X(-2,1,1,1,1)^{\oplus 2}\bigoplus \mathcal{O}_X(-1,0,0,0,0)^{\oplus 3}\arrow{d} & &  &\\
 & & \mathcal{O}_X^{\oplus 2} & & &
\end{tikzcd}
\end{center}

We have $\widetilde{B(-1,-2)} \cong \mathcal{I}/\mathcal{I}^2$. Hence, by Lemma \ref{l:gammadladywizorow}, $\Gamma_*(\mathcal{I}/\mathcal{I}^2) \cong B(-1,-2)$. 
Using Theorem \ref{t:pullback} and Lemma \ref{l:gammadladywizorow} we conclude that $\Gamma_*(G^*\mathcal{I}/\mathcal{I}^2)\cong R(-4,1,1,1,1)$. Therefore, from Lemma \ref{l:gladkieMDS_morfizmy_pochodza_z_algebry}, it follows that there is a graded homomorphism of $\mathrm{Cl}(X)$-graded $R$-modules $\eta: R(-4,1,1,1,1)\to \Gamma_*(G^*\Omega_{\mathbb{P}^1\times\mathbb{P}^2}^1)$ such that $\widetilde{\operatorname{coker}\eta}$ is isomorphic to $\Omega_X^1$. Moreover, $\eta$ is injective since it was obtained from an injective map of quasicoherent sheaves by applying first the pullback functor along a flat morphism and
then the global section functor. Denote the cokernel of $\eta$ by $P$. We have two exact sequences:
\[
0 \to R(-4,1,1,1,1)\xrightarrow{\eta} \Gamma_*(G^*\Omega_{\mathbb{P}^1\times\mathbb{P}^2}^1) \to P\to 0
\]
and
\begin{equation}\label{e:ciag_del_pezzo}
0 \to M_R \oplus N_R \to R^{\oplus 2}(-2,1,1,1,1)\oplus R^{\oplus 3}(-1,0,0,0,0) \xrightarrow{\epsilon\otimes_B 1_R} R^{\oplus 2}.
\end{equation}
The second sequence is exact on the left since, as can be checked, $\operatorname{Tor}_1^B(\operatorname{im{\epsilon}}, R)=0$. 
Note that by Lemma \ref{l:dla_gladkich_moge_wziac_dowolny} we have $\widetilde{M_R\oplus N_R}\cong G^*\Omega_{\mathbb{P}^1\times\mathbb{P}^2}^1$. Therefore, the same argument as before for $\Omega_{\mathbb{P}^1\times\mathbb{P}^2}^1$ shows, that $\Gamma_*(G^*\Omega_{\mathbb{P}^1\times\mathbb{P}^2}^1)$ is the kernel of  $\epsilon\otimes_B 1_R$. Thus, $\Gamma_*(G^*\Omega_{\mathbb{P}^1\times\mathbb{P}^2}^1)$ is isomorphic with $M_R\oplus N_R$ as a $\mathrm{Cl}(X)$-graded $R$-module.

Let $H(\alpha, \beta,x,y,z)=\alpha(xz-yz)+\beta(xy-yz)$ denote the equation of $Z$. Consider the graded homomorphisms of $\mathrm{Cl}(X)$-graded $R$-modules $\chi:R(-4,1,1,1,1)\to R^{\oplus 2}(-2,1,1,1,1)\oplus R^{\oplus 3}(-1,0,0,0,0)$ given by 
\[
f\mapsto ( \overline{G}^*(\frac{\partial H}{\partial \alpha}) f,\overline{G}^*(\frac{\partial H}{\partial \beta}) f,\overline{G}^*(\frac{\partial H}{\partial x}) f, \overline{G}^*(\frac{\partial H}{\partial y}) f,\overline{G}^*(\frac{\partial H}{\partial z}) f).
\]

This map is injective and its composition with $\epsilon\otimes_B 1_R$ is zero. Hence it factors through $\zeta:R(-4,1,1,1,1) \to M_R\oplus N_R$.
From the above considerations, it follows that $P\cong (M_R\oplus N_R) / \zeta(R)$. We will identify $M_R\oplus N_R$ with its image in $R^{\oplus 2}(-2,1,1,1,1)\oplus R^{\oplus 3}(-1,0,0,0,0)$ and denote the quotient $(M_R\oplus N_R)/\chi(R(-4,1,1,1,1))$ by $Q$.

From Lemma \ref{l:dla_gladkich_moge_wziac_dowolny}, it follows that $F^*\Omega_X^1 \cong \widetilde{Q \otimes_R A}$. 
After applying the functor $(.)\otimes_R A$ to sequence (\ref{e:ciag_del_pezzo}), we obtain an exact sequence of graded homomorphisms of $\mathbb{Z}$-graded $A$-modules:
\[
0\to M_A \oplus N_A \to A^{\oplus 2}\oplus A^{\oplus 3}(-2) \to A\oplus A
\]
where again the exactness on the left follows from the equality $\operatorname{Tor}^B_1(\operatorname{im}\epsilon, A)=0$. 

We can identify $M_A\oplus N_A$ inside $A^{\oplus 2}\oplus A^{\oplus 3}(-2)$ with the free $A$-module generated by $v_1=(0,0,2b, -a-b, 0)$, $v_2=(0,0,0,-a-b, 2a)$ and $v_3 = (1,2,0,0,0)$. This gives an identification of $M_A \oplus N_A$ with $A(-3)\oplus A(-3) \oplus A$.
Using this identification, 
the map $\chi\otimes_R 1_A: A(-4)\to A(-3)\oplus A(-3)\oplus A \subseteq  A^{\oplus 2} \oplus A^{\oplus 3}(-2)$ is given by the matrix:
\[
\begin{pmatrix} 
-8 b \\
-8 a\\
16ab(a^2-b^2)\\						
\end{pmatrix}.
\]
The $\mathbb{Z}$-graded $A$-module we are interested in, $Q\otimes_R A$, is isomorphic to the cokernel of this map. Consider the map $A(-3)\oplus A(-3)\oplus A \to A(-2)\oplus A$ given by $(f,g,h)\mapsto (af-bg, h+2a(a^2-b^2)f)$. The kernel of this homomorphism is precisely the image of $\chi \otimes_R 1_A$. Therefore, it gives an injective graded homomorphism of $\mathbb{Z}$-graded $A$-modules $Q\otimes_R 1_A \to A(-2)\oplus A$. The image is $A(-2)_{\geq 3}\oplus A$. This proves that $F^* \Omega_X^1 \cong \mathcal{O}_{\mathbb{P}^1}(-2)\oplus \mathcal{O}_{\mathbb{P}^1}$.

\begin{remark} We have started with a $\mathrm{Cl}(X)$-graded $R$-module $Q$ that describes the cotangent sheaf $\Omega_X^1$. The $\mathbb{Z}$-graded $A$-module that we have obtained $Q\otimes_R A$ 
describes the pullback of this quasicoherent sheaf but it does not satisfy $\Gamma_*(\widetilde{Q \otimes_R A}) \cong Q \otimes_R A$.
\end{remark}

\providecommand{\WileyBibTextsc}{}
\let\textsc\WileyBibTextsc
\providecommand{\othercit}{}
\providecommand{\jr}[1]{#1}
\providecommand{\etal}{~et~al.}

\end{document}